\documentclass[12pt,letterpaper]{amsart}
\usepackage{graphicx}
\usepackage[linktocpage=true,colorlinks,citecolor=blue,linkcolor=blue,urlcolor=blue]{hyperref}
\usepackage{amssymb}
\usepackage{cite}
\usepackage{amsmath}
\usepackage{latexsym}
\usepackage[colorinlistoftodos]{todonotes}
\usepackage{amscd}
\usepackage{tikz}
\usepackage{amsthm}
\usepackage{mathrsfs}
\usepackage{url}
\usepackage[utf8]{inputenc}
\usepackage[english]{babel}
\usepackage{amsfonts}
\usepackage{mathtools}
\usepackage{comment} 
\usepackage[autostyle]{csquotes}
\usepackage{thm-restate}
\usepackage{todonotes}

\numberwithin{equation}{section}
\begin{document}
\author{Bryce Kerr}
\address{BK: School of Science, University of New South Wales, Canberra, Australia
}
\email{bryce.kerr@unsw.edu.au}

\author{Oleksiy Klurman}
\address{OK: School of Mathematics, University of Bristol, United Kingdom 
}
\email{lklurman@gmail.com}
%
    
\newtheorem{theorem}{Theorem}
\newtheorem{lemma}[theorem]{Lemma}
\newtheorem{example}[theorem]{Example}
\newtheorem{algol}{Algorithm}
\newtheorem{corollary}[theorem]{Corollary}
\newtheorem{prop}[theorem]{Proposition}
\newtheorem{proposition}[theorem]{Proposition}
\newtheorem{problem}[theorem]{Problem}
\newtheorem{conj}[theorem]{Conjecture}
\newtheorem{cor}[theorem]{Corollary}

\theoremstyle{remark}
\newtheorem{definition}[theorem]{Definition}
\newtheorem{question}[theorem]{Question}
\newtheorem{remark}[theorem]{Remark}
\newtheorem*{acknowledgement}{Acknowledgements}

\numberwithin{equation}{section}
\numberwithin{theorem}{section}
\numberwithin{table}{section}
\numberwithin{figure}{section}
\numberwithin{lemma}{section}
\numberwithin{prop}{section}
\numberwithin{cor}{section}
\allowdisplaybreaks

\definecolor{olive}{rgb}{0.3, 0.4, .1}
\definecolor{dgreen}{rgb}{0.,0.5,0.}

\def\cc#1{\textcolor{red}{#1}} 

\def\Cde{C_{d,e}}
\def\tCde{\widetilde{C}_{d,e}}

\definecolor{dgreen}{rgb}{0.,0.6,0.}
\def\tgreen#1{\begin{color}{dgreen}{\it{#1}}\end{color}}
\def\tblue#1{\begin{color}{blue}{\it{#1}}\end{color}}
\def\tred#1{\begin{color}{red}#1\end{color}}
\def\tmagenta#1{\begin{color}{magenta}{\it{#1}}\end{color}}
\def\tNavyBlue#1{\begin{color}{NavyBlue}{\it{#1}}\end{color}}
\def\tMaroon#1{\begin{color}{Maroon}{\it{#1}}\end{color}}

\def\ndiv{\nmid}
\def \balpha{\bm{\alpha}}
\def \bbeta{\bm{\beta}}
\def \bgamma{\bm{\gamma}}
\def \bdelta{\bm{\delta}}
\def \blambda{\bm{\lambda}}
\def \bchi{\bm{\chi}}
\def \bphi{\bm{\varphi}}
\def \bpsi{\bm{\psi}}
\def \bnu{\bm{\nu}}
\def \bomega{\bm{\omega}}

\def\td{\widetilde d}
\def \te{\widetilde e} 
\def\talpha{\widetilde \alpha}
\def \tbeta{\widetilde \beta} 

\def \tcA{\widetilde {\cA}}
\def  \tcB{\widetilde{\cB}}

\def\vh{\vec{h}} 
\def\vj{\vec{j}} 

\def\vk{\vec{k}} 
\def\vl{\vec{l}} 

\def\vu{\vec{u}} 
\def\vv{\vec{v}} 

\def\vx{\vec{x}} 
\def\vy{\vec{y}} 
 
%


 \def\mand{\qquad\mbox{and}\qquad}

\def\cA{{\mathcal A}}
\def\cB{{\mathcal B}}
\def\cC{{\mathcal C}}
\def\cD{{\mathcal D}}
\def\cE{{\mathcal E}}
\def\cF{{\mathcal F}}
\def\cG{{\mathcal G}}
\def\cH{{\mathcal H}}
\def\cI{{\mathcal I}}
\def\cJ{{\mathcal J}}
\def\cK{{\mathcal K}}
\def\cL{{\mathcal L}}
\def\cM{{\mathcal M}}
\def\cN{{\mathcal N}}
\def\cO{{\mathcal O}}
\def\cP{{\mathcal P}}
\def\cQ{{\mathcal Q}}
\def\cR{{\mathcal R}}
\def\cS{{\mathcal S}}
\def\cT{{\mathcal T}}
\def\cU{{\mathcal U}}
\def\cV{{\mathcal V}}
\def\cW{{\mathcal W}}
\def\cX{{\mathcal X}}
\def\cY{{\mathcal Y}}
\def\cZ{{\mathcal Z}}

\def\E{\mathbb{E}}
\def\C{\mathbb{C}}
\def\F{\mathbb{F}}
\def\K{\mathbb{K}}
\def\Z{\mathbb{Z}}
\def\R{\mathbb{R}}
\def\Q{\mathbb{Q}}
\def\N{\mathbb{N}}
\def\L{\mathbb{L}}
\def\M{\textsf{M}}
\def\U{\mathbb{U}}
\def\P{\mathbb{P}}
\def\A{\mathbb{A}}
\def\fp{\mathfrak{p}}
\def\fq{\mathfrak{q}}
\def\n{\mathfrak{n}}
\def\X{\mathcal{X}}
\def\x{\textrm{\bf x}}
\def\w{\textrm{\bf w}}
\def\ovQ{\overline{\Q}}
\def \Kab{\K^{\mathrm{ab}}}
\def \Qab{\Q^{\mathrm{ab}}}
\def \Qtr{\Q^{\mathrm{tr}}}
\def \Kc{\K^{\mathrm{c}}}
\def \Qc{\Q^{\mathrm{c}}}
\def\ZK{\Z_\K}
\def\ZKS{\Z_{\K,\cS}}
\def\ZKSf{\Z_{\K,\cS_f}}
\def\RSf{R_{\cS_{f}}}
\def\RTf{R_{\cT_{f}}}

\def\S{\mathcal{S}}
\def\vec#1{\mathbf{#1}}
\def\ov#1{{\overline{#1}}}
\def\Sp{{\operatorname{S}}}
\def\Gm{\G_{\textup{m}}}
\def\fA{{\mathfrak A}}
\def\fB{{\mathfrak B}}
\def\fM{{\mathfrak M}}

\def \brho{\bm{\rho}}
\def \btau{\bm{\tau}}

\def\house#1{{%
    \setbox0=\hbox{$#1$}
    \vrule height \dimexpr\ht0+1.4pt width .5pt depth \dp0\relax
    \vrule height \dimexpr\ht0+1.4pt width \dimexpr\wd0+2pt depth \dimexpr-\ht0-1pt\relax
    \llap{$#1$\kern1pt}
    \vrule height \dimexpr\ht0+1.4pt width .5pt depth \dp0\relax}}


\newenvironment{notation}[0]{%
  \begin{list}%
    {}%
    {\setlength{\itemindent}{0pt}
     \setlength{\labelwidth}{1\parindent}
     \setlength{\labelsep}{\parindent}
     \setlength{\leftmargin}{2\parindent}
     \setlength{\itemsep}{0pt}
     }%
   }%
  {\end{list}}

\newenvironment{parts}[0]{%
  \begin{list}{}%
    {\setlength{\itemindent}{0pt}
     \setlength{\labelwidth}{1.5\parindent}
     \setlength{\labelsep}{.5\parindent}
     \setlength{\leftmargin}{2\parindent}
     \setlength{\itemsep}{0pt}
     }%
   }%
  {\end{list}}
\newcommand{\Part}[1]{\item[\upshape#1]}

\def\Case#1#2{%
\smallskip\paragraph{\textbf{\boldmath Case #1: #2.}}\hfil\break\ignorespaces}

\def\Subcase#1#2{%
\smallskip\paragraph{\textit{\boldmath Subcase #1: #2.}}\hfil\break\ignorespaces}

\renewcommand{\a}{\alpha}
\renewcommand{\b}{\beta}
\newcommand{\g}{\gamma}
\renewcommand{\d}{\delta}
\newcommand{\e}{\epsilon}
\newcommand{\f}{\varphi}
\newcommand{\fhat}{\hat\varphi}
\newcommand{\bfphi}{{\boldsymbol{\f}}}
\renewcommand{\l}{\lambda}
\renewcommand{\k}{\kappa}
\newcommand{\lhat}{\hat\lambda}
\newcommand{\bfmu}{{\boldsymbol{\mu}}}
\renewcommand{\o}{\omega}
\renewcommand{\r}{\rho}
\newcommand{\rbar}{{\ov\rho}}
\newcommand{\s}{\sigma}
\newcommand{\sbar}{{\ov\sigma}}
\renewcommand{\t}{\tau}
\newcommand{\z}{\zeta}


\newcommand{\ga}{{\mathfrak{a}}}
\newcommand{\gb}{{\mathfrak{b}}}
\newcommand{\gn}{{\mathfrak{n}}}
\newcommand{\gp}{{\mathfrak{p}}}
\newcommand{\gP}{{\mathfrak{P}}}
\newcommand{\gq}{{\mathfrak{q}}}

\newcommand{\Abar}{{\ov A}}
\newcommand{\Ebar}{{\ov E}}
\newcommand{\kbar}{{\ov k}}
\newcommand{\Kbar}{{\ov K}}
\newcommand{\Pbar}{{\ov P}}
\newcommand{\Sbar}{{\ov S}}
\newcommand{\Tbar}{{\ov T}}
\newcommand{\gbar}{{\ov\gamma}}
\newcommand{\lbar}{{\ov\lambda}}
\newcommand{\ybar}{{\ov y}}
\newcommand{\phibar}{{\ov\f}}

\newcommand{\Acal}{{\mathcal A}}
\newcommand{\Bcal}{{\mathcal B}}
\newcommand{\Ccal}{{\mathcal C}}
\newcommand{\Dcal}{{\mathcal D}}
\newcommand{\Ecal}{{\mathcal E}}
\newcommand{\Fcal}{{\mathcal F}}
\newcommand{\Gcal}{{\mathcal G}}
\newcommand{\Hcal}{{\mathcal H}}
\newcommand{\Ical}{{\mathcal I}}
\newcommand{\Jcal}{{\mathcal J}}
\newcommand{\Kcal}{{\mathcal K}}
\newcommand{\Lcal}{{\mathcal L}}
\newcommand{\Mcal}{{\mathcal M}}
\newcommand{\Ncal}{{\mathcal N}}
\newcommand{\Ocal}{{\mathcal O}}
\newcommand{\Pcal}{{\mathcal P}}
\newcommand{\Qcal}{{\mathcal Q}}
\newcommand{\Rcal}{{\mathcal R}}
\newcommand{\Scal}{{\mathcal S}}
\newcommand{\Tcal}{{\mathcal T}}
\newcommand{\Ucal}{{\mathcal U}}
\newcommand{\Vcal}{{\mathcal V}}
\newcommand{\Wcal}{{\mathcal W}}
\newcommand{\Xcal}{{\mathcal X}}
\newcommand{\Ycal}{{\mathcal Y}}
\newcommand{\Zcal}{{\mathcal Z}}

\renewcommand{\AA}{\mathbb{A}}
\newcommand{\BB}{\mathbb{B}}
\newcommand{\CC}{\mathbb{C}}
\newcommand{\FF}{\mathbb{F}}
\newcommand{\GG}{\mathbb{G}}
\newcommand{\KK}{\mathbb{K}}
\newcommand{\NN}{\mathbb{N}}
\newcommand{\PP}{\mathbb{P}}
\newcommand{\QQ}{\mathbb{Q}}
\newcommand{\RR}{\mathbb{R}}
\newcommand{\ZZ}{\mathbb{Z}}

\newcommand{\bfa}{{\boldsymbol a}}
\newcommand{\bfb}{{\boldsymbol b}}
\newcommand{\bfc}{{\boldsymbol c}}
\newcommand{\bfd}{{\boldsymbol d}}
\newcommand{\bfe}{{\boldsymbol e}}
\newcommand{\bff}{{\boldsymbol f}}
\newcommand{\bfg}{{\boldsymbol g}}
\newcommand{\bfi}{{\boldsymbol i}}
\newcommand{\bfj}{{\boldsymbol j}}
\newcommand{\bfk}{{\boldsymbol k}}
\newcommand{\bfm}{{\boldsymbol m}}
\newcommand{\bfp}{{\boldsymbol p}}
\newcommand{\bfr}{{\boldsymbol r}}
\newcommand{\bfs}{{\boldsymbol s}}
\newcommand{\bft}{{\boldsymbol t}}
\newcommand{\bfu}{{\boldsymbol u}}
\newcommand{\bfv}{{\boldsymbol v}}
\newcommand{\bfw}{{\boldsymbol w}}
\newcommand{\bfx}{{\boldsymbol x}}
\newcommand{\bfy}{{\boldsymbol y}}
\newcommand{\bfz}{{\boldsymbol z}}
\newcommand{\bfA}{{\boldsymbol A}}
\newcommand{\bfF}{{\boldsymbol F}}
\newcommand{\bfB}{{\boldsymbol B}}
\newcommand{\bfD}{{\boldsymbol D}}
\newcommand{\bfG}{{\boldsymbol G}}
\newcommand{\bfI}{{\boldsymbol I}}
\newcommand{\bfM}{{\boldsymbol M}}
\newcommand{\bfP}{{\boldsymbol P}}
\newcommand{\bfX}{{\boldsymbol X}}
\newcommand{\bfY}{{\boldsymbol Y}}
\newcommand{\bfzero}{{\boldsymbol{0}}}
\newcommand{\bfone}{{\boldsymbol{1}}}

\newcommand{\aff}{{\textup{aff}}}
\newcommand{\Aut}{\operatorname{Aut}}
\newcommand{\Berk}{{\textup{Berk}}}
\newcommand{\Birat}{\operatorname{Birat}}
\newcommand{\characteristic}{\operatorname{char}}
\newcommand{\codim}{\operatorname{codim}}
\newcommand{\Crit}{\operatorname{Crit}}
\newcommand{\critwt}{\operatorname{critwt}} 
\newcommand{\Cycle}{\operatorname{Cycles}}
\newcommand{\diag}{\operatorname{diag}}
\newcommand{\Disc}{\operatorname{Disc}}
\newcommand{\Div}{\operatorname{Div}}
\newcommand{\Dom}{\operatorname{Dom}}
\newcommand{\End}{\operatorname{End}}
\newcommand{\ExtOrbit}{\mathcal{EO}} 
\newcommand{\Fbar}{{\ov{F}}}
\newcommand{\Fix}{\operatorname{Fix}}
\newcommand{\FOD}{\operatorname{FOD}}
\newcommand{\FOM}{\operatorname{FOM}}
\newcommand{\Gal}{\operatorname{Gal}}
\newcommand{\genus}{\operatorname{genus}}
\newcommand{\GITQuot}{/\!/}
\newcommand{\GL}{\operatorname{GL}}
\newcommand{\GR}{\operatorname{\mathcal{G\!R}}}
\newcommand{\Hom}{\operatorname{Hom}}
\newcommand{\Index}{\operatorname{Index}}
\newcommand{\Image}{\operatorname{Image}}
\newcommand{\Isom}{\operatorname{Isom}}
\newcommand{\hhat}{{\hat h}}
\newcommand{\Ker}{{\operatorname{ker}}}
\newcommand{\Ksep}{K^{\textup{sep}}}  
\newcommand{\lcm}{{\operatorname{lcm}}}
\newcommand{\LCM}{{\operatorname{LCM}}}
\newcommand{\Lift}{\operatorname{Lift}}
\newcommand{\limstar}{\lim\nolimits^*}
\newcommand{\limstarn}{\lim_{\hidewidth n\to\infty\hidewidth}{\!}^*{\,}}
\newcommand{\llog}{\log\log}
\newcommand{\logplus}{\log^{\scriptscriptstyle+}}
\newcommand{\Mat}{\operatorname{Mat}}
\newcommand{\maxplus}{\operatornamewithlimits{\textup{max}^{\scriptscriptstyle+}}}
\newcommand{\MOD}[1]{~(\textup{mod}~#1)}
\newcommand{\Mor}{\operatorname{Mor}}
\newcommand{\Moduli}{\mathcal{M}}
\newcommand{\Norm}{{\operatorname{\mathsf{N}}}}
\newcommand{\notdivide}{\nmid}
\newcommand{\normalsubgroup}{\triangleleft}
\newcommand{\NS}{\operatorname{NS}}
\newcommand{\onto}{\twoheadrightarrow}
\newcommand{\ord}{\operatorname{ord}}
\newcommand{\Orbit}{\mathcal{O}}
\newcommand{\Per}{\operatorname{Per}}
\newcommand{\Perp}{\operatorname{Perp}}
\newcommand{\PrePer}{\operatorname{PrePer}}
\newcommand{\PGL}{\operatorname{PGL}}
\newcommand{\Pic}{\operatorname{Pic}}
\newcommand{\Prob}{\operatorname{Prob}}
\newcommand{\Proj}{\operatorname{Proj}}
\newcommand{\Qbar}{{\ov{\QQ}}}
\newcommand{\rank}{\operatorname{rank}}
\newcommand{\Rat}{\operatorname{Rat}}
\newcommand{\Res}{{\operatorname{Res}}}
\newcommand{\Resultant}{\operatorname{Res}}
\renewcommand{\setminus}{\smallsetminus}
\newcommand{\sgn}{\operatorname{sgn}}
\newcommand{\SL}{\operatorname{SL}}
\newcommand{\Span}{\operatorname{Span}}
\newcommand{\Spec}{\operatorname{Spec}}
\renewcommand{\ss}{{\textup{ss}}}
\newcommand{\stab}{{\textup{stab}}}
\newcommand{\Stab}{\operatorname{Stab}}
\newcommand{\Support}{\operatorname{Supp}}
\newcommand{\Sym}{\operatorname{Sym}}  
\newcommand{\tors}{{\textup{tor}}}
\newcommand{\Trace}{\operatorname{Trace}}
\newcommand{\trianglebin}{\mathbin{\triangle}} 
\newcommand{\tr}{{\textup{tr}}} 
\newcommand{\UHP}{{\mathfrak{h}}}    
\newcommand{\Wander}{\operatorname{Wander}}
\newcommand{\<}{\langle}
\renewcommand{\>}{\rangle}

\newcommand{\pmodintext}[1]{~\textup{(mod}~#1\textup{)}}
\newcommand{\ds}{\displaystyle}
\newcommand{\longhookrightarrow}{\lhook\joinrel\longrightarrow}
\newcommand{\longonto}{\relbar\joinrel\twoheadrightarrow}
\newcommand{\SmallMatrix}[1]{%
  \left(\begin{smallmatrix} #1 \end{smallmatrix}\right)}
  
  \def\({\left(}
\def\){\right)}
\def\fl#1{\left\lfloor#1\right\rfloor}
\def\rf#1{\left\lceil#1\right\rceil}


\title
{How negative can $\sum_{n\le x}\frac{f(n)}{n}$ be?}

\begin{abstract} 
Tur\'an observed that logarithmic partial sums $\sum_{n\le x}\frac{f(n)}{n}$ of completely multiplicative functions (in the particular case of the Liouville function $f(n)=\lambda(n)$) tend to be positive. We develop a general approach to prove two results aiming to explain this phenomena.\\
Firstly, we show that for every $\varepsilon>0$ there exists some $x_0\ge 1,$ such that for any completely multiplicative function $f$ satisfying $-1\le f(n)\le 1$, we have 
$$\sum_{n\le x}\frac{f(n)}{n}\ge -\frac{1}{(\log\log{x})^{1-\varepsilon}}, \quad x\ge x_0.$$ This improves a previous bound due to Granville and Soundararajan.
Secondly, we show that if $f$ is a typical (random) completely multiplicative function $f:\mathbb{N}\to \{-1,1\}$, the probability that $\sum_{n\le x}\frac{f(n)}{n}$ is negative for a given large $x,$ is $O(\exp(-\exp(\frac{\log x\cdot \log\log\log x}{C\log \log x}))).$ This improves on recent work of Angelo and Xu.

\end{abstract}

\maketitle

\section{Introduction}
An immediate consequence of the celebrated Dirichlet's class number formula is that $L(1,\chi)>0$ for any primitive real multiplicative character $\chi$ of modulus $q$. It thus follows that there exists some $x_0(q)$ which may depend on $q$ such that 
\begin{align}
\label{eq:s-zero}
\sum_{n\le x}\frac{\chi(n)}{n}\ge 0, \quad \text{for all} \quad x\ge x_0(q).
\end{align}
Establishing quantitative variants of~\eqref{eq:s-zero} is a fundamental problem related to the existence of putative Siegel zeros. Such problems have origins in the work of 
Tur\'{a}n~\cite{Turan}, who showed that if  partial sums of the Liouville function satisfy
\begin{align}
\label{eq:lio}
\sum_{n\le x}\frac{\lambda(n)}{n}\ge 0, \quad \text{for all} \quad n \ge 1,
\end{align}
then the Riemann hypothesis is true. Haselgrove~\cite{Has} proved that, in fact~\eqref{eq:lio} is false, with a rather amusing
 $x = 72, 185, 376, 951, 205$ being 
the smallest integer counterexample to \eqref{eq:lio} which was found in~\cite{BFM}. This naturally raises the question: why is this number so large?

On the other hand, as noted by Granville and Soundararajan~\cite{GS}, the above considerations combined with quadratic reciprocity imply there exist real Dirichlet $L$-functions with negative truncations
 \begin{align}
\sum_{n\le x}\frac{\chi(n)}{n}< 0, \quad \text{for some} \quad x\ge x_0,
\end{align}
 for any large $x_0\ge 1.$
 In particular, $x_0$ in~\eqref{eq:s-zero} cannot be taken to be uniform with respect to $q$ and one may ask how negative the sums in~\eqref{eq:s-zero} can get.

To facilitate our discussion, as in~\cite{GS}, we let $\cF,\cF_1$ and $\cF_0$ denote the set of multiplicative functions satisfying $-1\le f(n)\le 1, f(n)=\pm 1$ and $f(n)\in \{-1,0,1\}$ respectively. A quantitative form of the above  question is to establish lower bounds for the quantities 
\begin{align*}
\delta(x)=\min_{f\in \cF}\sum_{n\le x}\frac{f(n)}{n}, \quad \delta_i(x)=\min_{f\in \cF_i}\sum_{n\le x}\frac{f(n)}{n}, \quad i\in \{-1,0,1\}.
\end{align*}
We clearly have
\begin{align}
\label{eq:deltaequiv}
\delta(x)\le \delta_0(x)\le \delta_1(x),
\end{align}
and a quadratic reciprocity argument yields
\begin{align*}
\delta_0(x)=\min_{\substack{\chi \\ \text{quadratic character}}}\sum_{n\le x}\frac{\chi(n)}{n}.
\end{align*}
In~\cite{GS}, Granville and Soundararajan showed that
for sufficiently large $x$, we have 
\begin{align}
\label{GS:main}
\delta(x)\ge -\frac{1}{(\log\log{x})^{3/5}},
\end{align}
and 
\begin{align*}
\delta_1(x)\le -\frac{c}{\log{x}},
\end{align*}
and raised the question of improving these results.
In view of~\eqref{eq:deltaequiv}, the latter furnishes upper and lower bounds for $\delta,\delta_0,\delta_1,$ however understanding their asymptotic behaviour remains a mystery. Our first result gives the following improvement of \eqref{GS:main}.
\begin{theorem}
\label{thm:main}
For any $\varepsilon>0$ there exists $x_0$ such that if  $x\ge x_0,$ then
\begin{align*}
\delta(x)\ge -\frac{1}{(\log\log{x})^{1-\varepsilon}}.
\end{align*}
\end{theorem}
Another route to explain the positivity bias is to ask, given large $x,$ how likely it is for $\sum_{n\le x}\frac{f(n)}{n}$ to be negative for a typical (random) completely multiplicative function? To this end, we let $(f(p))_{p\ prime}$ to be a sequence
of independent random variables taking values $\pm 1$ with probability $1/2$ and put $f(n)=\prod_{p^k\vert| n}f^k(p)$ for all $n\ge 1.$ Sign changes of partial sums of random multiplicative functions have been subject of a series of recent investigations, see~\cite{AX},\cite{AHZ} and~\cite{K}.\\ Exploring a connection to the Euler product approximation, Angelo and Xu~\cite{AX}  showed that 
if $f$ is a random completely multiplicative function $f:\mathbb{N}\to \{\pm 1\}$, the probability that $\sum_{n\le x}\frac{f(n)}{n}$ is negative for a given large $x,$ is extremely small, that is $O(\exp(-\exp(\frac{\log x}{C\log \log x}))),$ for some absolute constant $C>0.$ Using a different approach we establish somewhat stronger result.
\begin{theorem}
\label{thm:main3}
Let $f:\mathbb{N}\to\{-1,1\}$ be a random completely multiplicative function. Then for any large $x,$ the probability that 
\begin{align}
\label{random:neg}
\sum_{n\le x}\frac{f(n)}{n}<0,
\end{align}
is $O(\exp(-\exp(\frac{\log x\cdot \log\log\log x}{C\log \log x}))).$
\end{theorem}
On the other hand, Granville and Soundararajan showed that for any large $x,$ there is a function $f:\mathbb{N}\to \{-1,1\},$ for which~\eqref{random:neg} holds (see~\cite{GS}). Thus the probability of the event~\eqref{random:neg} is at least $\frac{1}{2^{\pi(x)}}.$ It remains an interesting open problem to determine the exact order of the above quantity. In this direction, we note it does not seem difficult to modify their argument to show~\eqref{random:neg} holds with probability $\frac{1}{2^{c\pi(x)}},$ for some $c<1$.
\subsection{Ideas of the proofs.} Roughly speaking, both proofs start with a simple convolution identity (which is also used in~\cite{GS} but not in~\cite{AX}) and our improvements come from introducing bilinear structure to incorporate finer information on the distribution of $f$ on large primes $x^v\le p\le x$. More precisely, if $g(n)=\sum_{d\vert n}f(n)$ then
\[\sum_{n\le x}\frac{f(n)}{n}=\frac{1}{x}\sum_{n\le x}g(n)+\frac{1}{x}\sum_{n\le x}f(n)\left\{\frac{x}{n}\right\},\]
and crucially $g(n)\ge 0$ for all $n\ge 1.$ The remaining task is to produce satisfactory lower and upper bounds for the first and second terms respectively.
In the case of Theorem~\ref{thm:main}, one may proceed further and write as in Lemma~\ref{lem:gs},
\[\sum_{n\le x}\frac{f(n)}{n}=\frac{1}{x}\sum_{n\le x}g(n)+(1-\gamma)\frac{1}{x}\sum_{n\le x}f(n)+O\left(\frac{1}{\log^{1/5} x}\right).\]
To upper bound the second term we use Hal\'asz type estimates due to Hall and Tenenbaum (Lemma~\ref{lem:ht1}) and more recent bounds due to Granville, Harper and Soundararajan (Lemma~\ref{lem:ht}). The key novelty in our approach lies in treating the lower bounds for the first sum. Here, rather than using a classical bound due to Hildebrand (as in \cite{GS}), we rely on a recent additive combinatorics result of Matom\"{a}ki and Shao~\cite{MS}, which allows us to prove the main technical Proposition~\ref{lem:lb}.\\ 
In order to prove Theorem 1.2, we notice that due to positivity
\[\frac{1}{x}\sum_{n\le x}g(n)\ge \frac{1}{x}\sum_{p\le x} (1+f(p)),\]
and direct application of Bernstein's inequality produces a simple lower bound $\frac{1}{x}\sum_{n\le x}g(n)\gg \frac{1}{\log x}$ with an acceptable probability. We are thus left with estimating the probability of the event
\[\sum_{n\le x}f(n)\left\{\frac{x}{n}\right\}\ll -\frac{x}{\log x}.\]
We do so by a direct application of the majorant principle  and upper bounding very high moments of the unweighted partial sums $\mathbb{E}|\sum_{n\le x }f(n)|^{2q}.$ The key observation here (see  Proposition~\ref{thm:main2}) is that since we are interested in the event when the partial sums are extremely large (in absolute value), it is more beneficial to work ``close to the $1-$ line" rather than proceed, as in the work of Harper~\cite{Harp}, to ``$1/2-$ line"  (where the main interest is when partial sums have typical size around $\sqrt{x}$). This is accomplished by introducing bilinear structure and using bounds on the density of smooth numbers as well as Rankin's trick to treat contribution of large primes.\\
We believe that Proposition~\ref{lem:lb} and Proposition~\ref{thm:main2} are of independent interest and could be useful in future investigations.

\section{Preliminaries}

Here we collect several key facts, based on the argument of Granville and Soundararajan~\cite{GS}.
The following is~\cite[Proposition~3.1]{GS}.
\begin{lemma}
\label{lem:gs}
Let $f\in \cF$ and define 
$g(n)=\sum_{d|n}f(d).$
Then 
\begin{align*}
\sum_{n\le x}\frac{f(n)}{n}=\frac{1}{x}\sum_{n\le x}g(n)+(1-\gamma)\frac{1}{x}\sum_{n\le x}f(n)+O\left(\frac{1}{(\log{x})^{1/5}}\right),
\end{align*}
where $\gamma$ denotes Euler's constant.
\end{lemma}
We will use a version of Hal\'{a}sz's theorem in the form given by Granville, Harper and Soundararajan~\cite{HGS}.
\begin{lemma}
\label{lem:ht}
Let $f\in \cF$ and define $M=M(x)$ by 
\begin{align*}
\max_{|t|\ll \log{x}}\left|\frac{F(1+1/\log{x}+it)}{1+1/\log{x}+it}\right|:=e^{-M}(\log{x}),
\end{align*}
where 
$F(s)=\sum_{n=1}^{\infty}\frac{f(n)}{n^{s}}\cdot$
Then
\begin{align*}
\sum_{n\le x}f(n)\ll (1+M)e^{-M}x+\frac{(\log{x})^{o(1)}}{(\log{x})}\cdot
\end{align*}
\end{lemma}
A result of Hall and Tenenbaum~\cite{HT} allows one to remove the dependence on $t$ in Lemma~\ref{lem:ht} at the cost of a worse upper bound.
\begin{lemma}
\label{lem:ht1}
For any $f\in \cF$  we have 
$$\sum_{n\le x}f(n)\ll x\exp\left(-\kappa\sum_{p\le x}\frac{1-f(p)}{p}\right),$$
with $\kappa=0.32867$.
\end{lemma}
\section{Lower bounds for nonegative multiplicative functions}
We next show how ideas from~\cite{GKM,MS} may be used to obtain lower bounds for the function $g(n)$ defined in Lemma~\ref{lem:gs}. The main result of this section is Proposition~\ref{lem:lb}.  

The following is due to Matom\"{a}ki and Shao~\cite[Hypothesis P]{MS} which improves work of Granville, Koukoulopoulos and Matom\"{a}ki~\cite{GKM}. 
\begin{lemma}
\label{lem:ms}
Fix $\lambda\in (0,1)$. If $x$ is sufficiently large, $u,v$ satisfy 
$$1\le u \le v \le \frac{\log{x}}{1000\log\log{x}},$$
 and $\cP$ is a subset of the primes in $(x^{1/v},x^{1/u}]$ with 
$$\sum_{\substack{p\in \cP}}\frac{1}{p}\ge \frac{1+\lambda}{u},$$
then there exists an integer $k\in [u,v]$ such that  
\begin{align*}
\left|\left\{ (p_1,\dots,p_k)\in \cP^{k} \ : \ \frac{x}{2}\le p_1\dots p_k\le x \right\} \right|\ge \pi_v \frac{x}{v^{k}\log{x}},
\end{align*}
where $\pi_v$ is a constant with $\pi_v=v^{-o(v)}$ as $v\rightarrow \infty$. If $u$ is fixed and $v\ge 1000u^2/\lambda^2$, one can take $k\le e^{-1/u}v$.
\end{lemma}
Our next result uses arguments in the spirit of~\cite[pg. 218]{Hil}.
\begin{lemma}
\label{lem:log}
Let $g$ be a nonnegative multiplicative function satisfying 
\begin{align}
\label{eq:gass}
g(p^{j})\le j+1,
\end{align}
for each prime power $p^{j}$. For any
  real numbers $x,z\ge 1$, we have 
\begin{align*}
\sum_{\substack{n\le x \\ p|n \implies p\le z}}\frac{g(n)}{n}\gg \exp\left(\sum_{p\le z}\frac{g(p)}{p}\right)\left(1-C\exp\left(-\frac{\log{x}}{\log{z}}\right)\right),
\end{align*}
for some absolute constant $C$.
\end{lemma}
\begin{proof}
Let 
$\delta=\frac{1}{\log{z}},$
and consider 
\begin{align*}
\sum_{\substack{n\le x \\ p|n \implies p\le z}}\frac{g(n)}{n}&\ge \sum_{\substack{n\ge 1 \\ p|n \implies p\le z}}\frac{g(n)}{n}-\sum_{\substack{n> x \\ p|n \implies p\le z}}\frac{g(n)}{n}  \\
&\ge \prod_{p\le z}\left(1+\sum_{j\ge 1}\frac{g(p^{j})}{p^{j}} \right)-\frac{1}{x^{\delta}}\sum_{\substack{n\ge 1 \\ p|n \implies p\le z}}\frac{g(n)}{n^{1-\delta}}.
\end{align*}
By~\eqref{eq:gass} we get
\begin{align*}
\log\left(\sum_{\substack{n\ge 1 \\ p|n \implies p\le z}}\frac{g(n)}{n^{1-\delta}}\right)&=\sum_{p\le z}\frac{g(p)}{p}+O(1).
\end{align*}
This implies that for some absolute constant $C,$
\begin{align*}
\sum_{\substack{n\le x \\ p|n \implies p\le z}}\frac{g(n)}{n}\gg \exp\left(\sum_{p\le z}\frac{g(p)}{p}\right)\left(1-C\exp\left(-\frac{\log{x}}{\log{z}}\right)\right).
\end{align*}
and the result follows.
\end{proof}
We are now in a position to establish the main result of this section.
\begin{prop}
\label{lem:lb}
Let $\varepsilon>0$ be sufficiently small, $f\in \cF$ and define 
$g(n)=\sum_{d|n}f(d).$
For $0<\delta <1$, let $\cP_{\delta}$ denote the set 
$\cP_{\delta}=\{ p \ \ \text{prime} \ : \ f(p)\ge -\delta \ \},$
and suppose for some 
$$\frac{40000}{\varepsilon^2}\le v \le \frac{\log{x}}{1000\log\log{x}},$$
we have 
\begin{align}
\label{eq:123123}
\sum_{\substack{p\in \cP_{\delta} \\ x^{1/v}\le p \le x}}\frac{1}{p}\ge 1+\varepsilon.
\end{align}
Then 
\begin{align*}
\sum_{n\le x}g(n) \gg \varepsilon^2\left(\frac{(1-\delta)}{v}\right)^{v(1+o(1))/e}\exp\left(\sum_{p\le x}\frac{f(p)}{p}\right)x.
\end{align*} 
\end{prop}
\begin{proof}
 Let 
$G=\sum_{n\le x}g(n),$
and define
\begin{align}
\label{eq:ABdef}
\cA=\cP_{-1}\cap [1,x^{1/v}], \quad \cB=\cP_{\delta}\cap (x^{1/v},x],
\end{align}
so that 
\begin{align*}
G& \ge \sum_{\substack{a\le x^{\varepsilon/5} \\ p|a \implies p\in \cA}}g(a)\sum_{\substack{b\le x/a \\ p|b \implies p\in \cB}}g(b).
\end{align*}
If $p\in \cP_{\delta},$ then once again
$g(p^{\alpha})\ge 1-\delta,$
and hence 
\begin{align}
\label{eq:Glb1}
G\ge \sum_{\substack{a\le x^{\varepsilon/5} \\ p|a \implies p\in \cA}}g(a)\sum_{\substack{b\le x/a \\ p|b \implies p\in \cB}}(1-\delta)^{\Omega(b)},
\end{align}
where $\Omega(n)$ counts the number of prime factors of an integer $n$ with multiplicity.

 Our next step is to apply Lemma~\ref{lem:ms} to the inner summation over $b$. This requires verifying for each $a\le x^{\varepsilon/5}$ the lower bound 
\begin{align}
\label{eq:a}
\sum_{\substack{ p\in \cP_{\delta} \\ (x/a)^{v}\le p \le x/a}}\frac{1}{p}\ge 1+\frac{\varepsilon}{2}.
\end{align}
By~\eqref{eq:123123} 
\begin{align*}
\sum_{\substack{ p\in \cP_{\delta} \\ (x/a)^{v}\le p \le x/a}}\frac{1}{p}\ge \sum_{\substack{ p\in \cP_{\delta} \\ x^{v}\le p \le x}}\frac{1}{p}-\sum_{\substack{ p\in \cP_{\delta} \\ x/a\le p \le x }}\frac{1}{p}\ge 1+\varepsilon-\sum_{\substack{ p\in \cP_{\delta} \\ x/a\le p \le x }}\frac{1}{p},
\end{align*}
and for $a$ satisfying $a\le x^{\varepsilon/5}$ we have 
\begin{align*}
\sum_{\substack{ p\in \cP_{\delta} \\ x/a\le p \le x }}\frac{1}{p}&\le \log\log{x}-\log\log{(x/a)}+o(1) \\ & \le \log\log{x}-\log\log{x^{1-\varepsilon/5}}\le \frac{\varepsilon}{4}+o(1),
\end{align*}
which combined with the above implies~\eqref{eq:a}. Hence by Lemma~\ref{lem:ms}, for each such $a$, if $v$ satisfies 
\begin{align}
\label{eq:vconds1-1}
v\ge \frac{40000}{\varepsilon^2},
\end{align}
 there exists an integer 
$k\le \frac{v}{e},$
such that 
$$\left|\left\{ (p_1,\dots,p_k)\in \cP_{\delta}^{k} \ : \ \frac{x}{2a}\le p_1\dots p_k \le \frac{x}{a} \right\} \right| \ge v^{o(v)}\frac{x}{v^{k}a\log{x}}.$$
We see that
\begin{align*}
\sum_{\substack{b\le x/a \\ p|b \implies p\in \cB}}(1-\delta)^{\Omega(b)}&\ge v^{o(v)}\left(\frac{(1-\delta)}{v}\right)^{k}\frac{x}{a \log{x}} \\ & \ge \left(\frac{(1-\delta)}{v}\right)^{v(1+o(1))/e}\frac{x}{a \log{x}},
\end{align*}
and hence by~\eqref{eq:Glb1}
\begin{align}
\label{eq:G123}
G\ge \left(\frac{(1-\delta)}{v}\right)^{v(1+o(1))/e}\frac{x}{\log{x}}\sum_{\substack{a\le x^{\varepsilon/5} \\ p|a \implies p\in \cA}}\frac{g(a)}{a}.
\end{align}
Recalling~\eqref{eq:ABdef}
\begin{align*}
\sum_{\substack{a\le x^{\varepsilon/5} \\ p|a \implies p\in \cA}}\frac{g(a)}{a}=\sum_{\substack{a\le x^{\varepsilon/5} \\ p|a \implies p\le x^{1/v}}}\frac{g(a)}{a},
\end{align*}
which combined with~\eqref{eq:vconds1-1} and Proposition~\ref{lem:lb} yields
\begin{align*}
\sum_{\substack{a\le x^{\varepsilon/5} \\ p|a \implies p\in \cA}}\frac{g(a)}{a}\gg \exp\left(\sum_{p\le x^{\varepsilon^2/40000}}\frac{g(p)}{p}\right)\left(1+O(\varepsilon)\right).
\end{align*}
Assuming $\varepsilon$ sufficiently small and using~\eqref{eq:G123} gives the bound
\begin{align*}
G\gg \left(\frac{(1-\delta)}{v}\right)^{v(1+o(1))/e}\frac{x}{\log{x}}\exp\left(\sum_{p\le x^{\varepsilon^2/40000}}\frac{g(p)}{p}\right),
\end{align*}
and the result follows after noting 
\begin{align*}
\sum_{p\le x^{\varepsilon^2/40000}}\frac{g(p)}{p}&\ge \sum_{p\le x}\frac{g(p)}{p}+2\log{\varepsilon}+O(1) \\ 
&\ge \log\log{x}+\sum_{p\le x}\frac{f(p)}{p}+2\log{\varepsilon}+O(1).
\end{align*}
\end{proof}
\section{Proof of Theorem~\ref{thm:main}}

Since we may assume 
\begin{align*}
\sum_{n\le x}\frac{f(n)}{n}\le -\frac{1}{(\log{x})^{1/6}},
\end{align*}
 Lemma~\ref{lem:gs} and Lemma~\ref{lem:ht} imply  
\begin{align}
\label{eq:gs0}
\sum_{n\le x}\frac{f(n)}{n}\ge 2\left(\frac{1}{x}\sum_{n\le x}g(n)-C_0\left|\sum_{n\le x}f(n)\right|\right),
\end{align}
and
\begin{align}
\label{eq:gs1}
\sum_{n\le x}\frac{f(n)}{n}\ge 2\left(\frac{1}{x}\sum_{n\le x}g(n)-C(1+M)e^{-M}\right),
\end{align}
for  absolute constants $C_0,C$ and  
\begin{align}
\label{eq:M123}
\left|\frac{F(1+1/\log{x}+it)}{1+1/\log{x}+it}\right|=e^{-M}(\log{x}),
\end{align}
for some $|t|\ll \log{x}.$
 Let $\varepsilon,\varepsilon_1,\varepsilon_2>0$ be small and put
\begin{align}
\label{eq:deltav}
\delta=1-\varepsilon_1, \quad v=(\log\log{x})^{1-\varepsilon_2},
\end{align}
and consider the set
\begin{align*}
\cP_{\delta}=\{p\le x \ : \ f(p)\ge -\delta\}.
\end{align*}
We distinguish between
\begin{align}
\label{eq:case1}
\sum_{\substack{x^{1/v}\le p \le x \\ p\in \cP_{\delta}}}\frac{1}{p}\le  1+\varepsilon,
\end{align}
and
\begin{align}
\label{eq:case2}
\sum_{\substack{x^{1/v}\le p \le x \\ p\in \cP_{\delta}}}\frac{1}{p}> 1+\varepsilon.
\end{align}
If~\eqref{eq:case1} holds, then 
\begin{align*}
\Re\left(\sum_{p\le x}\frac{1-f(p)p^{-it}}{p}\right)&\ge \Re\left(\sum_{\substack{x^{1/v}\le p \le x \\ p\not \in \cP_{\delta}}}\frac{1-f(p)p^{-it}}{p}\right)+O(1) \\ 
& \ge \Re\left(\sum_{\substack{x^{1/v}\le p \le x \\ p\not \in \cP_{\delta}}}\frac{1+p^{-it}}{p}\right)-\sum_{\substack{x^{1/v}\le p \le x \\ p\not \in \cP_{\delta}}}\frac{|1+f(p)|}{p}+O(1).
\end{align*}
If $p\not \in \cP_{\delta},$ then 
$|1+f(p)|\le 1-\delta \le \varepsilon_1,$
and consequently 
\begin{align*}
\Re\left(\sum_{p\le x}\frac{1-f(p)p^{-it}}{p}\right)&\ge \Re\left(\sum_{\substack{x^{1/v}\le p \le x }}\frac{1+p^{-it}}{p}\right)-\varepsilon_1 \log{v}+O(1) \\
&\ge (1-\varepsilon_1)\log{v}+O(1).
\end{align*}
The above implies that 
\begin{align*}
\frac{1}{\log{x}}\left|\frac{F(1+1/\log{x}+it)}{1+1/\log{x}+it}\right|\ll \frac{1}{v^{1-\varepsilon_1}},
\end{align*}
which combined with~\eqref{eq:gs1},~\eqref{eq:M123} and~\eqref{eq:deltav} gives
\begin{align}
\label{eq:case11}
\sum_{n\le x}\frac{f(n)}{n}\ge -\frac{C}{(\log\log{x})^{(1-\varepsilon_1)(1-\varepsilon_2)}}.
\end{align}
If~\eqref{eq:case2} holds, then by~\eqref{eq:gs1} and Proposition~\ref{lem:lb}
\begin{align*}
\sum_{n\le x}g(n)\gg \varepsilon^2\left(\frac{(1-\delta)}{v}\right)^{v(1+o(1))/e}\log{x}\exp\left(-\sum_{p\le x}\frac{1-f(p)}{p}\right) .
\end{align*}
Recalling~\eqref{eq:deltav}, we have 
\begin{align*}
\left(\frac{(1-\delta)}{v}\right)^{v(1+o(1))/e}\log{x}&\ge \exp\left(\log\log{x}-v\log{v/\varepsilon_1}\right) \\
&\gg \exp(\log\log{x}-(\log\log{x})^{1-\varepsilon_1/2}\gg (\log{x})^{1-o(1)}.
\end{align*}
Upon applying Lemma~\ref{lem:ht1} and~\eqref{eq:gs0} we deduce
\begin{align*}
\sum_{n\le x}\frac{f(n)}{n}\ge C_1(\log{x})^{1-o(1)}\exp\left(-\sum_{p\le x}\frac{1-f(p)}{p}\right)-C_2\exp\left(-\kappa\sum_{p\le x}\frac{1-f(p)}{p}\right),
\end{align*}
for some absolute constants $C_1,C_2$. This implies
\begin{align*}
\sum_{n\le x}\frac{f(n)}{n}\ge -\frac{c}{\log\log{x}},
\end{align*}
and combined with~\eqref{eq:case11} we complete the proof, after taking $\varepsilon_1,\varepsilon_2$ sufficiently small and $x$ sufficiently large.
\section{Random multiplicative functions}

Recall that $\cF_1$ denotes the space of completely multiplicative functions $f$ satisfying $f(p)=\pm 1$ for $p\le x.$ Let $\mu$ denote the counting measure on $\cF_1$.

\begin{cor}
\label{thm:main1}
For any $M\gg 1$, there exists a constant $C$ such that 
\begin{align*}
\mu\left(\left\{ f\in \cF_1 \ : \ \left|\sum_{n\le N}\left\{\frac{x}{n}\right\}f(n)\right|\ge \frac{M}{\log{x}}x\right\} \right)\le \exp \left(-\exp\left(\frac{(\log \log \log{x})\log{x}}{C\log\log{x}}\right)\right).
\end{align*}
\end{cor}
Corollary~\ref{thm:main1} is a consequence of the following moment inequality.
\begin{prop}
\label{thm:main2}

 For any $y\le x,$  even integer $q$ and constant $C$, we have 
\begin{align*}
\E\left[\left(\sum_{n\le x}f(n)\right)^{q} \right]&\ll \frac{x^q}{y^{q/2}}\exp\left(q\left(q\exp\left(-\frac{(\log\log\log{x})\log{x}}{C\log\log{x}}\right)+2\log\log{x}\right)\right) \\ &+\left(\frac{x}{(\log{x})^{C}}\right)^{q}.
\end{align*}

\end{prop}
{\it Deduction of Theorem~\ref{thm:main3} from Corollary~\ref{thm:main1}.} Recall that, if $g(n)=\sum_{d\vert n}f(n)$ then
\[\sum_{n\le x}\frac{f(n)}{n}=\frac{1}{x}\sum_{n\le x}g(n)+\frac{1}{x}\sum_{n\le x}f(n)\left\{\frac{x}{n}\right\}.\]
Upon applying  Bernstein inequality, we have that 
\[\sum_{p\le x}f(p)\le -\frac{x}{2\log x}\]
with probability $O(\exp(-\frac{x}{100\log x})).$
Since $g(n)\ge 0$ for all $n\ge 1,$ we have that for large $x,$ 
\[\sum_{n\le x}g(n)\ge \frac{1}{x}\sum_{p\le x} (1+f(p))\ge \frac{1}{5\log x}\]
with an acceptable probability $1-O(\exp(-\frac{x}{100\log x})).$ We are left to apply Corollary~\ref{thm:main1} with $M=1/10,$ say and the result follows.

{\it Deduction of Corollary~\ref{thm:main1} from Proposition~\ref{thm:main2}.} By the majorant principle, for any even integer $q$
\begin{align*}
\E\left[\left(\sum_{n\le x}\left\{\frac{x}{n}\right\}f(n)\right)^{q} \right]\le \E\left[\left(\sum_{n\le x}f(n)\right)^{q} \right].
\end{align*}
Now, Corollary~\ref{thm:main1} follows from Proposition~\ref{thm:main2} by taking $y=(\log{x})^{A}$ and $q=\exp\left(\frac{(\log\log\log{x})(\log{x})}{C\log\log{x}}\right)$ and applying Markov's inequality.
\section{Proof of Proposition~\ref{thm:main2}} 
For $y\le x$ define
\begin{align}
\label{eq:delta}
\delta=\frac{\log{y}}{\log{x}},\quad \varepsilon=\frac{\log\log\log{x}}{C\log\log{x}}\cdot
\end{align}
Partition summation over $n$ as 
\begin{align*}
\sum_{\substack{n\le x \\ p|n \implies p>x^{\varepsilon}}}f(n)=\sum_{\substack{n\ell\le x \\ p|n \implies p\ge x^{\varepsilon} \\ p|\ell \implies  p < x^{\varepsilon}}}f(\ell)f(n)=\sum_{j\le \log{x}+1}\sum_{\substack{n\ell\le x \\ p|n \implies p\ge x^{\varepsilon} \\ p|\ell \implies p < x^{\varepsilon} \\ e^{j}\le \ell < e^{j+1}}}f(\ell)f(n),
\end{align*}
and apply Minkowski's inequality and the majorant principle to obtain
\begin{align*}
 \E\left[\left(\sum_{n\le x}f(n)\right)^{q} \right]^{1/q}\le \sum_{j\le \log{x}+1}\E\left[\left(\sum_{\substack{n\le x/e^{j} \\ p|n \implies p\ge x^{\varepsilon}}}f(n)\sum_{\substack{\ell \le e^{j+1} \\ p|\ell \implies  p\le  x^{\varepsilon}}}f(\ell)\right)^{q} \right]^{1/q}.
\end{align*}
Partition summation over $j$ into ``small" and ``large" in a natural fashion:
\begin{align}
\label{eq:ES12}
 \E\left[\left(\sum_{n\le x}f(n)\right)^{q} \right]^{1/q}\le S_1+S_2,
\end{align}
with 
\begin{align}
\label{eq:S1}
S_1=\sum_{j\le \log{x}/2}\E\left[\left(\sum_{\substack{n\le x/e^{j} \\ p|n \implies p\ge x^{\varepsilon}}}f(n)\sum_{\substack{\ell \le e^{j+1} \\ p|\ell \implies  p<  x^{\varepsilon}}}f(\ell)\right)^{q} \right]^{1/q}
\end{align}
\begin{align*}
S_2=\sum_{\log{x}/2<j\le \log{x}+1}\E\left[\left(\sum_{\substack{n\le x/e^{j} \\ p|n \implies p\ge x^{\varepsilon}}}f(n)\sum_{\substack{\ell \le e^{j+1} \\ p|\ell \implies  p<  x^{\varepsilon}}}f(\ell)\right)^{q} \right]^{1/q}.
\end{align*}
Consider first $S_2$. Our plan is to bound inner summation trivially by using density of smooth numbers to establish a satisfactory bound. Indeed, there exists $j\ge \log{x}/2$ such that 
\begin{align*}
S_2\le (\log{x})\E\left[\left(\sum_{\substack{n\le x/e^{j} \\ p|n \implies p\ge x^{\varepsilon}}}f(n)\sum_{\substack{\ell \le e^{j+1} \\ p|\ell \implies  p<  x^{\varepsilon}}}f(\ell)\right)^{q} \right]^{1/q}.
\end{align*} 
Define $\beta$ by $x^{\beta}=e^{j+1}$, so that 
$\beta\ge 1/2.$
 We have 
\begin{align*}
\sum_{\substack{n\le x/e^{j} \\ p|n \implies p\ge x^{\varepsilon}}}f(n)\sum_{\substack{\ell \le e^{j+1} \\ p|\ell \implies  p<  x^{\varepsilon}}}f(\ell)\ll x^{1-\beta}\Psi(x^{\beta},x^{\varepsilon}),
\end{align*}
where 
$$\Psi(x,y)=\sum_{\substack{n\le x \\ p|n \implies p\le y}}1.
$$
This implies  
\begin{align}
\label{eq:S2}
S_2\ll (\log{x})x^{1-\beta}\Psi(x^{\beta},x^{\varepsilon}).
\end{align}
By a result of Canfield, Erd\"{o}s and Pomerance~\cite{CEP}, for $u=\frac{\log{w}}{\log{z}},$ we have
\begin{align}
\label{eq:smooth1}
\psi(w,z)\le \frac{w}{u^{1+o(1)}}, \quad \text{provided} \quad u\le z^{1-o(1)} \quad \text{and} \quad u\rightarrow \infty. 
\end{align}
We apply~\eqref{eq:smooth1} with 
$
w=x^{\beta}, z=x^{\varepsilon}.
$
For this choice of parameters
$u=\frac{\beta}{\varepsilon},$
and we note that by~\eqref{eq:delta} $u\rightarrow \infty$ and 
\begin{align*}
u=\frac{\beta}{\varepsilon}\le \frac{C\log\log{x}}{\log\log\log{x}}\le x^{C\log\log\log{x}/\log\log{x}}\le x^{\varepsilon/2}.
\end{align*}
Combining~\eqref{eq:S2} and~\eqref{eq:smooth1} we get
\begin{align*}
S_2&\ll (\log{x})\left(\frac{\varepsilon}{\beta}\right)^{\beta/\varepsilon}x \ll  \frac{x}{(\log{x})^{C/4}}.
\end{align*}
Consequently, using~\eqref{eq:ES12}
\begin{align}
\label{eq:S12}
 \E\left[\left(\sum_{n\le x}f(n)\right)^{q} \right]^{1/q}\ll S_1+\frac{x}{(\log{x})^{C/4}}.
\end{align}
We now turn our attention to estimating $S_1.$ With $\delta,\epsilon$ as in~\eqref{eq:delta}, another application of the majorant principle gives 

\begin{align*}
S_1\ll \frac{x}{y} \sum_{j\le \log{x}/2}y^{j/(\log{x})}\E\left[\left(\sum_{\substack{n\le x/e^{j} \\ p|n \implies p\ge x^{\varepsilon}}}\frac{f(n)}{n^{1-\delta}}\sum_{\substack{\ell \le e^{j+1} \\ p|\ell \implies  p< x^{\varepsilon}}}\frac{f(\ell)}{\ell}\right)^{q} \right]^{1/q},
\end{align*}
and hence by applying majorant principle once again we arrive at
\begin{align}
\label{eq:S1}
\nonumber S_1&\ll \frac{x}{y^{1/2}}\E\left[\left(\sum_{\substack{ p|n \implies x^{\varepsilon}\le p\le x}}\frac{f(n)}{n^{1-\delta}}\sum_{\substack{ p|\ell \implies  p<  x^{\varepsilon}}}\frac{f(\ell)}{\ell}\right)^{q} \right]^{1/q} \\ 
&=\frac{x}{y^{1/2}}\E\left[\left(\prod_{x^{\varepsilon}\le p\le x}\left(1-\frac{f(p)}{p^{1-\delta}}\right)^{-1}\right)^{q} \right]^{1/q}\E\left[\left(\prod_{p\le x^{\varepsilon}}\left(1-\frac{f(p)}{p}\right)^{-1}\right)^{q} \right]^{1/q}.
\end{align}
We next expand expectation as a product,
\begin{align*}
\E\left[\left(\prod_{x^{\varepsilon}\le p\le x}\left(1-\frac{f(p)}{p^{1-\delta}}\right)^{-1}\right)^{q} \right]&\le \exp(O(q))\E\left[\left(\prod_{x^{\varepsilon}\le p\le x}\left(1+\frac{f(p)}{p^{1-\delta}}\right)\right)^{q} \right] .
\end{align*}
Observe that
\begin{align*}
\E\left[\left(\prod_{p\le x}\left(1+\frac{f(p)}{p^{1-\delta}}\right)\right)^{q} \right]=\prod_{p\le x}\frac{1}{2}\left( \left(1+\frac{1}{p^{1-\delta}}\right)^{q}+\left(1-\frac{1}{p^{1-\delta}}\right)^{q} \right).
\end{align*}

\begin{align*}
\prod_{x^{\varepsilon}\le p\le x}\frac{1}{2}\left( \left(1+\frac{1}{p^{1-\delta}}\right)^{q}+\left(1-\frac{1}{p^{1-\delta}}\right)^{q} \right)&\le \prod_{x^{\varepsilon}\le p \le x}\frac{\exp(q/p^{1-\delta})+\exp(-q/p^{1-\delta})}{2} \\ 
&\le \exp\left(\frac{q^2}{2x^{(1-2\delta)\varepsilon}}\right),
\end{align*}
and therefore
\begin{align*}
\E\left[\left(\prod_{x^{\varepsilon}\le p\le x}\left(1-\frac{f(p)}{p^{1-\delta}}\right)^{-1}\right)^{q} \right]&\le \exp\left(\frac{q^2}{2x^{\varepsilon/2}}+O(q)\right).
\end{align*}
After recalling~\eqref{eq:delta}, the latter yields
\begin{align}
\label{eq:E2}
\E\left[\left(\prod_{x^{\varepsilon}\le p\le x}\left(1-\frac{f(p)}{p}\right)^{-1}\right)^{q} \right]&\le \exp\left(q^2\exp\left(-\frac{(\log\log\log{x})(\log{x})}{4C\log\log{x}}\right)+O(q)\right).
\end{align}
We finally deduce the bound
\begin{align*}
\E\left[\left(\prod_{p\le x^{\varepsilon}}\left(1-\frac{f(p)}{p}\right)^{-1}\right)^{q} \right]&\ll \exp\left(O(q)+q\sum_{p\le x^{\varepsilon}}\frac{1}{p}\right) \\ &
\le \exp\left(O(q)+q\log \log{x}\right).
\end{align*}
Combining the above with~\eqref{eq:S1} and~\eqref{eq:E2}, we get 
\begin{align*}
S_1\ll \frac{x}{y^{1/2}}\exp\left(\left(q\exp\left(-\frac{(\log\log\log{x})(\log{x})}{4C\log\log{x}}\right)+2\log\log{x}\right)\right),
\end{align*}
and the result follows from~\eqref{eq:S12} after renaming the constant $C$.
\section{Acknowledgement}
The first author is currently supported by the Australian Research Council DE220100859 and part of this work was carried out while both authors were visiting the Max Planck Institute for Mathematics, Bonn. We thank MPIM (Bonn) for providing excellent working conditions and Focused Research Grant (HIMR) for the support.

\end{document}